
\magnification = 1200
\topskip = 8mm


\font \hugebf = cmbx10 at 14.4 pt
\font \hugebfit = cmbxti10 at 14.4pt
\font \bigbf = cmbx10 at 12pt
\font \bigbfit = cmbxti10 at 12pt
\font \bigrm = cmr10 at 12pt

\catcode`\@=11

\font\tenmsx=msam10
\font\sevenmsx=msam7
\font\fivemsx=msam5
\font\tenmsy=msbm10
\font\sevenmsy=msbm7
\font\fivemsy=msbm5
\newfam\msxfam
\newfam\msyfam
\textfont\msxfam=\tenmsx  \scriptfont\msxfam=\sevenmsx
 \scriptscriptfont\msxfam=\fivemsx
\textfont\msyfam=\tenmsy  \scriptfont\msyfam=\sevenmsy
 \scriptscriptfont\msyfam=\fivemsy

\font\tenmsa=msam10
\font\sevenmsa=msam7
\font\fivemsa=msam5
\newfam\msafam
  \textfont\msafam=\tenmsa
  \scriptfont\msafam=\sevenmsa
  \scriptscriptfont\msafam=\fivemsa
\def\msa{\fam\msafam\tenmsa}

\font\tenBbb=msbm10
\font\sevenBbb=msbm7
\font\fiveBbb=msbm5
\newfam\Bbbfam
  \textfont\Bbbfam=\tenBbb
  \scriptfont\Bbbfam=\sevenBbb
  \scriptscriptfont\Bbbfam=\fiveBbb

\def\hexnumber@#1{\ifnum#1<10 \number#1\else
 \ifnum#1=10 A\else\ifnum#1=11 B\else\ifnum#1=12 C\else
 \ifnum#1=13 D\else\ifnum#1=14 E\else\ifnum#1=15 F\fi\fi\fi\fi\fi\fi\fi}

\def\msx@{\hexnumber@\msxfam}
\def\msy@{\hexnumber@\msyfam}
\mathchardef\nmid="3\msy@2D
\mathchardef\varnothing="0\msy@3F
\mathchardef\nexists="0\msy@40
\mathchardef\smallsetminus="2\msy@72
\def\Bbb{\ifmmode\let\next\Bbb@\else
\def\next{\errmessage{Use \string\Bbb\space only in math mode}}\fi\next}
\def\Bbb@#1{{\Bbb@@{#1}}}
\def\Bbb@@#1{\fam\msyfam#1}

\catcode`\@=12

\def \bC{{\Bbb C}}
\def \bQ{{\Bbb Q}}

\def\Ker{\mathop{\rm Ker}\nolimits}

\def \noi {\noindent}
\def \dis {\displaystyle}

\def \nmid {\not\>|\>}
\def\compact{\mathop{\msa b}\nolimits}
\def\ssm{\mathop{\Bbb r}\nolimits}
\def\hugeLtwo{\hbox{\hugebfit L}^{\hbox{\bf 2}}}
\def\bigLtwo{\hbox{\bigbfit L}^{\hbox{\bf 2}}}

\def \cqfd{\unskip\kern 6pt\penalty 500
\raise -2pt\hbox{\vrule\vbox to10pt{\hrule width 4pt
\vfill\hrule}\vrule}\par}

\def\adots{\mathinner{\mkern2mu\raise1pt\hbox{.}
\mkern3mu\raise4pt\hbox{.}\mkern1mu\raise7pt\hbox{.}}}

\def\build#1_#2^#3{\mathrel{\mathop{\kern 0pt#1}\limits_{#2}^{#3}}}
\def\hfl#1#2{\smash{\mathop{\hbox to
8mm{\rightarrowfill}}\limits^{\scriptstyle#1}_{
\scriptstyle#2}}}
\def\vfl#1#2{\llap{$\scriptstyle #1$}\left\downarrow \vbox to
4mm{}\right.\rlap{$\scriptstyle #2$}}

\headline={\ifnum\pageno=1 {\hfill} \else{\hss \tenrm -- \folio\ -- \hss} \fi}
\footline={\hfil}

\catcode`\@=11
\newcount\@tempcnta \newcount\@tempcntb 
\def\heure{{%
\@tempcnta=\time \divide\@tempcnta by 60 \@tempcntb=\@tempcnta
\multiply\@tempcntb by -60 \advance\@tempcntb by \time
\ifnum\@tempcntb > 9 \number\@tempcnta:\number\@tempcntb
  \else\number\@tempcnta:0\number\@tempcntb\fi}}
\catcode`\@=12

\def\ermini{\raise 1ex\hbox{\pc{}er|}}
\def\aujourdhui{\ifnum\day=1 1\ermini\else\number\day\fi\
\ifcase\month\or janvier\or f\'evrier\or mars\or avril\or mai\or juin\or
juillet\or aout\or septembre\or octobre\or novembre\or d\'ecembre\fi\
\number\year}

\centerline{\hugebf Cohomologie $\hugeLtwo$ sur les rev\^etements}
\medskip
\centerline{\hugebf d'une vari\'et\'e complexe compacte}
\vskip 1cm
\centerline{\bigbf Fr\'ed\'eric Campana}\medskip
\centerline{\bigrm Universit\'e de Nancy I, D\'epartement de Math\'ematiques}
\bigskip
\centerline{\bigbf Jean-Pierre Demailly}\medskip
\centerline{\bigrm Universit\'e de Grenoble I, Institut Fourier}
\vskip 1cm
\centerline{Version du 27 janvier 2000, imprim\'ee le \aujourdhui, \heure}
\vskip 2cm

\noi {\bigbf Introduction.}  \bigskip \noi 
Andreotti-Vesentini [AV], Ohsawa [Oh], Gromov [G], Koll\'ar [K], entre
autres, ont montr\'e que la th\'eorie de Hodge d'une vari\'et\'e
k\"ahl\'erienne compacte pouvait \^etre d\'efinie avec les m\^emes
propri\'et\'es dans le cadre $L^2$ si cette vari\'et\'e \'etait
seulement compl\`ete. Par ailleurs, les th\'eor\`emes d'annulation de
la g\'eom\'etrie k\"ahl\'erienne ou projective reposant sur la
m\'ethode de Kodaira-Bochner-Nakano admettent par nature des versions
$L^2$ (voir Androtti-Vesentini [AV] et [D]).  On se propose ici de
d\'efinir une cohomologie $L^2$ naturelle sur tout rev\^etement
\'etale d'un espace analytique complexe $X$, \`a valeurs dans le
rel\`evement de tout faisceau analytique coh\'erent ${\cal F}$ sur
$X$.  Cette cohomologie a toutes les propri\'et\'es habituelles de la
cohomologie des faisceaux sur $X$ (suites exactes de cohomologie,
suites spectrales, th\'eor\`emes d'annulation, en particulier), et ces
propri\'et\'es sont obtenues en incorporant l'information issue des
estim\'ees $L^2$ dans les preuves standards des r\'esultats
correspondants.  La cohomologie $L^2$ devrait offrir un cadre
agr\'eable pour \'etudier la g\'eom\'etrie des rev\^etements, en
fournissant un formalisme fonctoriel jouissant des propri\'et\'es
attendues.  Lorsque l'espace $X$ de base est compact et que le
rev\^etement est galoisien de groupe $\Gamma$, on peut d\'efinir la
$\Gamma$-dimension des groupes de cohomologie $L^2$ associ\'es \`a un
faisceau coh\'erent sur la base. On \'etablit en particulier leur
finitude et on \'etend le th\'eor\`eme de l'indice $L^2$ de Atiyah
dans ce cadre. Enfin, si $X$ est projective, on a des th\'eor\`emes
d'annulation $L^2$ qui \'etendent naturellement les th\'eor\`emes
d'annulation usuels (th\'eor\`eme de Kodaira-Serre, th\'eor\`eme de
Kawamata-Viehweg~$\ldots$).
\medskip

Dans [E], P.~Eyssidieux a annonc\'e la construction d'une telle
cohomologie, en utilisant des proc\'ed\'es voisins de ceux pr\'esent\'es ici.
Notons aussi qu'un th\'eor\`eme d'annulation $L^2$ en cohomologie
$L^2$ similaire \`a 4.1 est \'enonc\'e par J.~Koll\'ar dans [K], 11.4.
\vfill\eject

\noi {\bigbf \S\ 1. Norme $\bigLtwo$ sur les sections.}

\bigskip

\noi {\bf 1.0.} Soit $X$ une vari\'et\'e analytique complexe, ${\cal F}$ un
faisceau analytique coh\'erent sur $X$, et $U$ un ouvert relativement compact
de $X$. On dira que $U$ est ${\cal F}${\it-admissible} s'il existe un ouvert
de {\it Stein} $V$ contenant $U$ et tel que $U$ soit relativement compact
dans $V$, ainsi qu'un morphisme surjectif $f : {\cal O}^r_V \to {\cal
F}_{|V}$ de faisceaux de ${\cal O}_V$-modules sur $V$. Un tel morphisme
sera appel\'e une $0$-pr\'esentation de~${\cal F}$. Si le fibr\'e vectoriel
trivial $V \times {\Bbb C}^r$ est muni d'une m\'etrique hermitienne $h$,
on d\'efinira, pour $s \in H^0(U,{\cal O}^r_U)$ la norme $\|s\|$ par :
$\|s\|^2 = \dis\int_U h(s,s) d\mu$, o\`u $\mu$ est la forme volume
d'une m\'etrique fix\'ee sur $X$. On notera $H^0_{(2)}(U,{\cal O}^r)$
l'espace vectoriel (de Hilbert) des $s$ tels que $\|s\| < +\infty$, et
$H^0_{(2)}(U,{\cal F}) := f_\ast H^0_{(2)}(U,{\cal O}^r) \subset
H^0(U,{\cal F})$.

\noi On notera que cet espace est ind\'ependant des choix $(h, \mu$,
et m\^eme $f$ -- voir ci-dessous) faits. On le munit de la norme $L^2$
quotient : pour $\sigma = f_\ast(s) \in H^0_{(2)}(U,{\cal F})$, on pose
$$\|\sigma\| := \inf \{\|s\| \mid f_\ast(s) = : f \circ s = \sigma,~ 
s \in H^0_{(2)}(U,{\cal O}^r)\}.$$

\bigskip

\proclaim {\bf 1.1.} Proposition. Si $\|\sigma\| = 0$, alors $\sigma =
0$.  Autrement dit, la semi-norme ainsi d\'efinie sur
$H^0_{(2)}(U,{\cal F})$ est une norme. De plus $H^0_{(2)}(U,{\cal F})$
\'equip\'e de cette norme est un espace de Hilbert isom\'etrique \`a
l'orthogonal $(\hbox{\rm Ker} f_\ast)^\bot$ de $\hbox{\rm Ker} f_\ast
: H^0_{(2)}(U,{\cal O}^r) \to H^0_{(2)}(U,{\cal F})$ dans
$H^0_{(2)}(U,{\cal O}^r)$, et $(\hbox{\rm Ker} f_\ast)$ est ferm\'e
dans $H^0_{(2)}(U,{\cal O}^r)$.

\bigskip

\noi {\bf D\'emonstration.} Il suffit de montrer que $\Ker f_\ast$ est
ferm\'e dans $H^0_{(2)}(U,{\cal O}^r)$. Or ceci r\'esulte du fait 
que la topologie $L^2$ est plus forte que la topologie de la convergence
uniforme sur les compacts de $U$ ([W], III.7), et du fait bien connu que le noyau
$\Ker f_\ast:H^0(U,{\cal O}^r) \to H^0(U,{\cal F})$ est ferm\'e pour
la topologie de la convergence uniforme sur les compacts (c'est le cas
pour les sections \`a valeurs dans un sous-faisceau quelconque,
[H], 6.3.5 et chap.~7).
\bigskip

\proclaim {\bf 1.2.} D\'efinition. Deux espaces de Hilbert $(E,h_i)$
$(i=1,2)$ sur le m\^eme espace sous-jacent $E$ sont dits \hbox{\it
\'equivalents} si les normes $h_1$ et $h_2$ d\'efinissent la m\^eme
topologie (ou encore : s'il existe $0 < A < B$ tels que : $Ah_1 < h_2
< Bh_1$).

\bigskip

\proclaim {\bf 1.3.} Corollaire. A \'equivalence pr\`es, l'espace de Hilbert
$\big(H^0_{(2)}(U,{\cal F}),\|\;\|\big)$ est ind\'ependant des choix
$(h\,;\mu\,;f)$ faits.

\bigskip

\noi {\bf D\'emonstration.} Seule l'ind\'ependance vis \`a vis de $f$ m\'erite
d'\^etre v\'erifi\'ee : soient $f_i : {\cal O}^{r_i}_V \!\to \!{\cal
F}_{|V} \!\to \!0$ $(i=1,2)$ deux $0$-pr\'esentations de ${\cal F}_{|V}$
sur un $V$ commun. Puisque $V$ est Stein, il existe $\varphi : {\cal
O}_V^{r_i}\to {\cal O}_V^{r_j}$, avec $i\not= j$, tel que $f_j \circ
\varphi = f_i$.

\noi Cette application fournit la continuit\'e de l'application identique de
$H^0_{(2)}(U,{\cal F})$ muni de la norme d\'eduite de $f_i$ dans lui-m\^eme
muni de la norme
d\'eduite de $f_j$.

\bigskip

\noi {\bf 1.4.} Soit $U' \subset U$ ; l'application naturelle de
restriction : $\hbox{res} :
H^0_{(2)}(U,{\cal F}) \to H^0_{(2)}(U',{\cal F})$ est continue, et compacte
si $U'\compact U$. L'affirmation est en effet claire dans le cas
de faisceaux localement libres, et le cas g\'en\'eral s'en d\'eduit.

\bigskip

\noi {\bf 1.5.} Soit $u:{\cal F}\to{\cal G}$ un morphisme de
faisceaux, et $U$ un ouvert qui soit \`a la fois ${\cal F}$-admissible
et ${\cal G}$-admissible relativement \`a un m\^eme ouvert de Stein
$V$.  On a alors un morphisme induit $u_{(2)}:H^0_{(2)}(U,{\cal F}) 
\to H^0_{(2)}(U,{\cal G})$ {\it continu}. En effet,
soit $f:{\cal O}_V^r\to{\cal F}_{|V}$ une $0$-pr\'esentation de ${\cal
F}_{|V}$ et $g:{\cal O}_V^{r+s}\to{\cal G}_{|V}$ une 
$0$-pr\'esentation de ${\cal G}_{|V}$ choisie en sorte que $g\circ
i=u\circ f$, o\`u $i:{\cal O}_V^r\to {\cal O}_V^{r+s}$ est l'injection
des $r$-premi\`eres composantes. Alors le noyau du morphisme
$f_\star:H_{(2)}^0(U,{\cal O}_V^r) \to H^0(U,{\cal F})$ s'envoie
par $i$ dans le noyau de $g_\star:H_{(2)}^0(U,{\cal
O}_V^{r+s})\to H^0(U,{\cal G})$, et on en d\'eduit le morphisme
$u_{(2)}$ voulu par passage au quotient. De plus $u_{(2)}$ est injectif si 
$u$ est injectif, et $u_{(2)}$ est surjectif si $u$ est surjectif.
\bigskip

\noi {\bf 1.6.} Le foncteur ${\cal F}\mapsto H_{(2)}^0(U,{\cal F})$ n'est
en g\'en\'eral pas exact. Pour le voir, on peut consid\'erer par exemple
le morphisme injectif $u:{\cal O}_{\bC^2}\to{\cal O}_{\bC^2}$,
$s\mapsto z_2s$. Alors le morphisme induit
$$u_{(2)}:H_{(2)}^0(U,{\cal O})\to H_{(2)}^0(U,{\cal O})$$
n'est pas d'image ferm\'ee sur la boule unit\'e $U=B(0,1)\subset\bC^2$
(on peut v\'erifier que la section $(1-z_1)^{-3/2}$ n'est pas dans $L^2(U)$,
tandis que $z_2(1-z_1)^{-3/2}$ est dans $L^2(U)$, et par suite
$z_2(1-z_1)^{-3/2}$ est seulement dans l'adh\'erence de l'image).
Ceci montre qu'on ne peut pas avoir une suite exacte
$$
H_{(2)}^0(U,{\cal O}_{\bC^2})\build\longrightarrow_{}^{\dis u_{(2)}}
H_{(2)}^0(U,{\cal O}_{\bC^2})\to H_{(2)}^0(U,{\cal O}_{\bC\times\{0\}}).
$$
\smallskip

\noi {\bf 1.7.} Le d\'efaut d'exactitude du foncteur sections $L^2$
pourra \^etre palli\'e par l'observation suivante: soit 
$$
{\cal F}\build\longrightarrow_{}^{u}{\cal G}
\build\longrightarrow_{}^{v}{\cal H}
$$
une suite exacte de faisceaux admettant des $0$-pr\'esentations sur
un ouvert de Stein $V$, et soient $U'\compact U\compact V$ des ouverts de
Stein. Il existe une constante $C>0$ telle que pour tout \'el\'ement
$g$ dans le noyau de $v_{(2)}:H^0_{(2)}(U,{\cal G})\longrightarrow H^0_{(2)}
(U,{\cal H})$, on puisse trouver un \'el\'ement $f\in H^0_{(2)}(U',{\cal F})$
tel que $u_{(2)}(f)=g_{|U'}$ et $\Vert f\Vert_{L^2(U')}\le C
\Vert g\Vert_{L^2(U)}$. En effet, la topologie de $H^0_{(2)}(U,{\cal G})$ 
est plus forte que la topologie de la convergence uniforme sur 
les compacts de $U$ (induite par passage au quotient \`a partir
d'une pr\'esentation ${\cal O}^N\to{\cal G}$ et de la topologie d'espace
de Fr\'echet sur $H^0(U,{\cal O}^N)$). On conclut \`a partir 
de la suite exacte d'espaces de Fr\'echet 
$$
H^0(U,{\cal F})\to H^0(U,{\cal G})\to H^0(U,{\cal H})
$$
et du fait que le morphisme de restriction $H^0(U,{\cal F})\to
H^0_{(2)}(U',{\cal F})$ est continu.
\bigskip

\noi {\bf 1.8. Remarque.} Les notions d'espaces de sections $L^2$ 
peuvent \'egalement \^etre d\'efinies de mani\`ere analogue pour un
espace analytique $X$ arbitraire (r\'eduit ou non), en plongeant localement
l'ouvert de Stein $V\subset X$ dans un espace ambiant lisse ${\Bbb C}^N$, 
et en consid\'erant l'extension triviale ${\cal G}$ du faisceau 
${\cal F}_{|V}$ \`a ${\Bbb C}^N$ (telle que ${\cal G}_{|{\Bbb C}^N
\setminus V}=0\,$). Les normes $L^2$ pour un ouvert $U\compact V$ sont 
alors calcul\'ees en travaillant sur un ouvert de Stein $U'\compact{\Bbb C}^N$ 
tel que $U'\cap V=U$.
\vskip0.75cm

\noi {\bigbf \S\ 2. Image directe $\bigLtwo$.}

\bigskip

\noi {\bf 2.0. Conventions.} Soit $X$ un espace analytique
complexe et $p : \widetilde X \to X$ un rev\^etement \'etale de
$X$. Soit ${\cal F}$ un faisceau analytique coh\'erent sur $X$ et
$\widetilde{\cal F} := p^\ast {\cal F}$ son rel\`evement \`a
$\widetilde X$.  Si $U \compact V$ sont des ouverts de $X$ avec
$V$ Stein et $f : {\cal O}^r_V \to {\cal F}_V$ une 0-r\'esolution de
${\cal F}$ sur $V$, on notera : $\widetilde U =: p^{-1}(U) \, ;
\widetilde V := p^{-1}(V)$ ; et $\tilde f : {\cal O}^r_{\widetilde V}
\to \widetilde{\cal F}_{\widetilde V}$ les rel\`evements
correspondants \`a $\widetilde X$.

\noi On dira que $V$ est $p${\it-simple} si chaque composante connexe
de $\widetilde V$ est appliqu\'ee par $p$ sur une composante connexe
de $V$, bijectivement. On supposera cette condition satisfaite.

\noi Soit $h$ une m\'etrique hermitienne sur le fibr\'e trivial $V
\times {\Bbb C}^r$ associ\'e \`a ${\cal O}^r_V$, et $\tilde h$ son
rel\`evement \`a $\widetilde V \times {\Bbb C}^r$, associ\'e \`a
${\cal O}^r_{\widetilde V}$.

\noi Ceci permet de d\'efinir la notion de norme $L^2$ pour $\tilde s
\in H^0(\widetilde U,\widetilde{\cal F})$, gr\^ace \`a la d\'efinition
de 1.0, et aussi $H^0_{(2)}(\widetilde U,\widetilde{\cal F})$ qui, 
muni de cette norme, est un espace de Hilbert.  De plus, les arguments
du {\S} 1 montrent que l'espace de Hilbert $\big(H^0_{(2)}(\widetilde
U,\widetilde{\cal F}),\|\;\|\big)$ est ind\'ependant des choix
$(V,f,h,\mu)$ faits, \`a \'equivalence pr\`es.

\bigskip

\proclaim {\bf 2.1.} D\'efinition. Soit $W$ un ouvert de $X$, et $\widetilde
W =: p^{-1}(W)$.
Soit $\tilde s \in H^0(\widetilde W,\widetilde{\cal F})$. On dit que
$\tilde s$ est \hbox{\it
localement $L^2$ sur $X$} si, pour chaque $x \in W$, il existe des
voisinage ouverts $U
\compact V$ de $x$ dans $X$, avec $V$ Stein et $p$-simple, tels que
la restriction de
$\tilde s$ \`a $\widetilde U$ soit dans $H^0_{(2)}(\widetilde
U,\widetilde{\cal F})$.

\bigskip

\noi {\bf 2.2.} L'ensemble, not\'e : $H^0_{(2),{\rm loc}}(\widetilde
W,\widetilde{\cal F})$ des $\tilde s$ de $H^0(\widetilde W,\widetilde{\cal
F})$ qui sont
localement $L^2$ sur $X$ forme clairement un espace vectoriel complexe. On
a de plus,
des applications naturelles de restriction res$_{W,W'} :
H^0_{(2),{\rm loc}}
(\widetilde W,\widetilde{\cal F}) \to H^0_{(2)}(\widetilde
W',\widetilde{\cal F})$ pour $W
\supset W'$, ouverts de $X$, et donc un pr\'efaisceau \`a valeurs dans la
cat\'egorie des
espaces vectoriels complexes : 
$$
W \to H^0_{(2),{\rm loc}}(\widetilde W,\widetilde{\cal F}).
$$ 
Il est imm\'ediat de v\'erifier que ce pr\'efaisceau est un faisceau
d'espaces vectoriels complexes sur $X$.
\bigskip

\proclaim {\bf 2.3.} D\'efinition. Le faisceau ainsi d\'efini : $W \to
H^0_{(2),{\rm loc}}(\widetilde W,\widetilde{\cal F})$ sur $X$ est not\'e
$p_{\ast(2)}\widetilde{\cal F}$ ; il
est appel\'e \hbox{\it le faisceau image directe} $L^2$ de $\widetilde{\cal
F}$ par $p$.

\bigskip

\noi {\bf 2.4.} On munit maintenant naturellement
$\big(p_{\ast(2)}\widetilde{\cal F}\big)$ d'une structure de ${\cal
O}_X$-module comme suit : si $\tilde s \in H^0_{(2),{\rm
loc}}(\widetilde W,\widetilde{\cal F})$ repr\'esente un germe de
section de $\big(p_{\ast(2)}\widetilde{\cal F}\big)$ en $x \in W$, et
si $\varphi \in H^0(W,{\cal O}_W)$, alors $(p^\ast\varphi\cdot\tilde
s) \in H^0_{(2),{\rm loc}}(\widetilde W,\widetilde{\cal F})$, qui est 
donc un $H^0(W,{\cal O}_W)$-module.

\bigskip

\noi {\bf 2.5. Remarque.} Le faisceau $p_{\ast(2)}(\widetilde{\cal F})$
n'est en g\'en\'eral coh\'erent que si $p$ est fini, auquel cas
il se r\'eduit \`a l'image directe $p_\ast\widetilde{\cal F}$ usuelle.
Si $p$ est infini et ${\cal F} \not= 0$, alors 
$p_{\ast(2)}(\widetilde{\cal F})$ n'est jamais coh\'erent.
\bigskip

\proclaim {\bf 2.6.} Proposition. Soit ${\cal F} \to {\cal G} \to {\cal H}$ une
suite exacte de faisceaux coh\'erents analytiques sur $X$ ; alors la suite
naturelle d'images
directes $L^2 : p_{\ast(2)}\widetilde{\cal F} \to
p_{\ast(2)}\widetilde{\cal G} \to
p_{\ast(2)}\widetilde{\cal H}$ est exacte. $($Autrement dit : le
foncteur d'image directe $L^2$ par $p$ est exact$)$.

\bigskip

\noi {\bf D\'emonstration.} Soit $V$ un ouvert $p$-simple sur lequel
${\cal F}$, ${\cal G}$ et ${\cal H}$ admettent des
$0$-pr\'e\-sentations.  et $U'\compact U \compact V$ des ouverts de
Stein connexes comme dans 1.7. Toutes les composantes
connexes $U'_j\compact U_j\compact V_j$ de $\widetilde U'$,
$\widetilde U$, $\widetilde V$ sont alors en isomorphisme avec
$U'\compact U\compact V$. Il existe par cons\'equent une constante $C$
ind\'ependante de $j$ telles que les sections $g_j$ du noyau de
$H^0_{(2)}(U_j,\widetilde{\cal G})\to H^0_{(2)}(U_j,\widetilde{\cal
H})$ se rel\`event en des sections $f_j$ de $H^0_{(2)}(U'_j,{\cal F})$,
avec $\Vert f_j\Vert_{L^2(U'_j)}\le C\Vert g_j
\Vert_{L^2(U_j)}$. Toute section $g=\bigoplus
g_j$ dans le noyau de $H^0_{(2)}(\widetilde U,\widetilde{\cal G})\to
H^0_{(2)}(\widetilde U,\widetilde {\cal H})$ se rel\`eve donc en une section
$f=\bigoplus f_j\in H^0_{(2)}(\widetilde U',\widetilde{\cal F})$
dans~$L^2(\widetilde U')$.
\bigskip

\proclaim {\bf 2.7.} Corollaire. On a un isomorphisme naturel de faisceaux 
de ${\cal O}_X$-modules :\break $\eta :
p_{\ast(2)}\widetilde{\cal F} \buildrel \sim \over {\to}
\big(p_{\ast(2)}{\cal O}_{\widetilde X}\big) \dis\otimes_{{\cal
O}_X} {\cal F}$.

\bigskip

\noi {\bf D\'emonstration.} C'est imm\'ediat \`a partir des
d\'efinitions lorsque ${\cal F}$ est localement libre. En
g\'en\'eral, on conclut \`a partir de l'exactitude du foncteur image
directe $L^2$, par r\'ecurrence sur la longueur d'une r\'esolution
libre (locale) de ${\cal F}$ sur $X$.

\bigskip

\noi {\bf 2.8. Exemple :} Supposons que ${\cal F}$ admette une r\'esolution
finie localement
libre (si $X$ est projective, c'est toujours le cas, par un r\'esultat de
J.-P. Serre) : $0 \to
{\cal L}_n \to {\cal L}_{n-1} \to \cdots \to {\cal L}_0 \to {\cal F} \to 0$.
Alors :
$p_{\ast(2)}\widetilde{\cal F}$ admet une r\'esolution finie de m\^eme longueur
par des
faisceaux localement de la forme $\big(p_{\ast(2)}{\cal O}_{\widetilde
X}\big)^{\oplus r}$
:
$$0 \to p_{\ast(2)}\widetilde{\cal L}_n \to \cdots \to
p_{\ast(2)}\widetilde{\cal L}_0 \to
p_{\ast(2)}{\cal F} \to 0.$$
Lorsque $X$ est projective, les ${\cal L}_k$ peuvent \^etre pris de la forme :
$${\cal L} = \bigoplus^r_{j=1} {\cal O}_X(a_j), \quad a_j \in {\Bbb Z}$$
et on a donc :
$$p_{\ast(2)}{\cal L} = \bigoplus^r_{j=1} p_{\ast(2)}{\cal O}_{\widetilde X}
\otimes_{{\cal O}_X} {\cal O}_X(a_j).$$
\bigskip

\proclaim {\bf 2.9.} Proposition. Soit $Y$ un sous-espace de
$X$, ${\cal F}'$ un faisceau analytique coh\'erent sur $Y$, et ${\cal
F}=i_\ast{\cal F}'$ son extension par $0$ sur $X\ssm Y$ (image directe
par l'injection $i:Y\to X$). Alors, si $p' : \widetilde Y:=p^{-1}(Y)
\to Y$ est le rev\^etement \'etale de $Y$ induit par $p : \widetilde X
\to X$, on a $p_{\ast(2)}{\widetilde{\cal F}}=\widetilde i_\ast(
p'_{\ast(2)}{\widetilde{\cal F}'})$ o\`u 
$\widetilde i:\widetilde Y\to\widetilde X$ est l'injection naturelle.
\bigskip

\noi {\bf D\'emonstration.} On observe d'abord que si $V$ est un
ouvert de Stein simple dans $X$ et si $U\compact V$ est un voisinage
dans $X$ d'un point $x\in Y$ (resp.\ $U'\compact Y\cap U$ un voisinage
de $x$ dans $Y$), on a un morphisme de restriction continu
$H^0_{(2)}(U,{\cal F})\to H^0_{(2)} (U',{\cal F'})$. En consid\'erant
les sections sur les rev\^etements $\widetilde U$ et $\widetilde U'$
puis en passant \`a la limite inductive sur $U$ et $U'$, on en
d\'eduit qu'on a un morphisme de restriction $p_{\ast(2)}{\cal F}\to
\widetilde i_\ast(p'_{\ast(2)}{\cal F}')$ continu.  Il s'agit en fait
d'un isomorphisme. En effet, prenons des voisinages ouverts de Stein
$U$, $U_1$ de $x$ dans $X$ (resp.\ $U'$ de $x$ dans $Y$), tels que
$U\compact U_1\compact V$ et $Y\cap U_1\compact U'$. On a un
homomorphisme surjectif d'espaces de Fr\'echet $H^0(U_1,{\cal F})\to
H^0(Y\cap U_1,{\cal F'})$, qui est par suite un morphisme ouvert. Pour
toute section $s'$ de ${\cal F'}$ sur $U'$, on peut alors trouver une
section $s$ de ${\cal F}$ sur $U_1$ telle que $s_{|Y\cap
U_1}=s'_{|Y\cap U_1}$ et $\Vert s\Vert_{L^\infty(U)}\le C\Vert
s'\Vert_{L^\infty(U'_1)}$ pour un certain $U'_1\compact Y\cap
U_1$. Au niveau des normes $L^2$, ceci implique $\Vert
s\Vert_{L^2(U)}\le C \Vert s'\Vert_{L^2(U')}$.  On conclut en passant
aux sections sur les rev\^etements $\widetilde U$ et $\widetilde U'$.
\vskip 0.75cm

\noi {\bigbf \S 3. Cohomologie $\bigLtwo$.}

\bigskip

\noi {\bf 3.0.} Les conventions sont celles de 2.0.

\bigskip

\proclaim {\bf 3.1.} D\'efinition. Soit $W$ un ouvert de $X$. On
d\'efinit la cohomologie $L^2$ de $\widetilde W$ \`a valeurs dans
$\widetilde{\cal F}$ comme \'etant celle du faisceau
$\big(p_{\ast(2)}\widetilde{\cal F}\big)$ sur $W$. On note alors
$H^q_{(2)}(\widetilde W,\widetilde{\cal F})$ $(q \ge 0)$ le $q$-i\`eme
groupe de cohomologie ainsi d\'efini.

\bigskip

\proclaim {\bf 3.2.} Th\'eor\`eme. Soit $0 \to {\cal F} \to {\cal G} \to {\cal
H} \to 0$ une suite
exacte de faisceaux analytiques coh\'erents sur $X$. On peut lui associer une suite
exacte longue de
cohomologie $L^2$, fonctorielle en ${\cal F}$
$$0 \to H^0_{(2)}(\widetilde X,\widetilde{\cal F}) \to \cdots \to
H^q_{(2)}(\widetilde X,\widetilde{\cal F}) \to
H^q_{(2)}(\widetilde X,\widetilde{\cal G}) \to H^q_{(2)}(\widetilde X,\widetilde{\cal H}) \to
H^{q+1}_{(2)}(\widetilde X,\widetilde{\cal F}) \to \cdots$$
\bigskip

\noi {\bf D\'emonstration.} Cette suite exacte  r\'esulte imm\'ediatement de l'exactitude du foncteur
$p_{\ast(2)}$ et des propri\'et\'es usuelles de la cohomologie.
\bigskip

\proclaim{\bf 3.3.} Proposition. Soit un diagramme
commutatif
$$
\matrix{
\widetilde Y & \buildrel {\displaystyle\widetilde f}\over 
\longrightarrow & \widetilde X\cr
\noalign{\vskip7pt}
\kern-10pt p'\Big\downarrow& &\Big\downarrow p\kern-5pt\cr
Y & \buildrel {\displaystyle f}\over \longrightarrow & X\cr
}
$$
o\`u les fl\`eches verticales sont des rev\^etements, $f:Y\to X$ un morphisme
analytique et $p':\widetilde Y\to Y$ le rev\^etement image inverse de $p$
par~$f$. Soit ${\cal F}$ un faisceau analytique coh\'erent sur $Y$. Alors 
pour tout $q\ge0$ on a la formule de commutation
$$
p_{\ast (2)}(R^qf_\ast{\cal F})^{\widetilde{}} =
R^q\widetilde f_\ast(p'_{\ast (2)}{\widetilde{\cal F}}).
$$
\bigskip

\noi {\bf D\'emonstration.} Les deux faisceaux en question sont les
faisceaux associ\'es au pr\'efaisceau
$$
U \longmapsto H^q(f^{-1}(U), p'_{\ast (2)}{\widetilde{\cal F}}),\qquad
U \longmapsto H^q_{(2)}(\widetilde f^{-1}(\widetilde U), \widetilde{\cal F}),
$$
et ces deux pr\'efaisceaux co\"{\i}ncident par d\'efinition de la cohomologie 
$L^2$. \hbox{\cqfd} 
\bigskip

\noi {\bf 3.4. Isomorphisme de Dolbeault.} On suppose que $X$ est une
vari\'et\'e lisse, que ${\cal F}$ est localement libre, et on note $F$ le 
fibr\'e vectoriel (suppos\'e muni d'une 
m\'etrique hermitienne) associ\'e sur $X$. Soit $\widetilde{\cal F}^{r,q}$
le faisceau des formes diff\'erentielles $v$ de type $(r,q)$ \`a
valeurs dans $\widetilde F$ et \`a
coefficients $L^2_{{\rm loc}}$ sur $\widetilde X$, telles que
$\overline\partial v$ soient aussi $L^2_{{\rm loc}}$.
Soit $p_{\ast(2)}\widetilde
{\cal F}^{r,q}$ le (pr\'e)faisceau sur $X$ d\'efini par $U \to
H^0_{(2),{\rm loc}}(\widetilde U,\widetilde{\cal F}^{r,q})$, 
image directe $L^2$ de 
$\widetilde{\cal F}^{r,q}$ sur~$X$, \`a savoir le faisceau des formes 
diff\'erentielles qui sont $L^2$ localement au dessus de $X$, sur
les ouverts de la forme~$\widetilde U=p^{-1}(U)$. L'op\'erateur $\overline
\partial$ fournit un complexe de faisceaux sur $X$
$$0 \to p_{\ast(2)}\widetilde{\cal F} \to p_{\ast(2)}\widetilde{\cal
F}^{0,0} \dis \mathop
{\longrightarrow}^{\overline\partial} \cdots \dis \mathop
{\longrightarrow}^{\overline\partial} p_{\ast(2)}\widetilde{\cal F}^{0,q} \to
p_{\ast(2)}\widetilde{\cal F}^{0,q+1} \to \cdots \to
p_{\ast(2)}\widetilde{\cal F}^{0,n} \to
0.$$
Ce complexe est exact, en vertu du th\'eor\`eme d'existence de 
H\"ormander-Andreotti-Vesen\-tini pour les solutions $L^2$ de l'op\'erateur
$\overline\partial$: en effet, on va appliquer ce th\'eor\`eme sur des ouverts 
$\widetilde U$ qui rev\^etent un ouvert de Stein $U\subset X$ quelconque,
en prenant une m\'etrique k\"ahl\'erienne $\omega$ sur $X$ et des poids
de la forme $\widetilde\varphi=\varphi\circ p$, o\`u
$\varphi$ est choisi en sorte que $i\partial\overline\partial\varphi+
{\rm Ricci}_\omega+{\rm Courbure}_F\ge \omega$. Pour toute forme $w$ telle
que $\overline \partial w = 0$, on voit alors que
l'\'equation $\overline\partial v = w$ admet une solution $v$ telle
que $\|v\| \le \|w\|$ en norme $L^2$.

On obtient ainsi une r\'esolution fine de $p_{\ast(2)}\widetilde{\cal F}$
(par des faisceaux de modules sur le faisceau d'anneaux des fonctions
$C^\infty$ sur $X$), de sorte que le complexe de Dolbeault $L^2$ calcule
bien la cohomologie $L^2$ d\'efinie en~3.1. Nous pouvons \'enoncer:

\proclaim {\bf 3.5.} Proposition. {\rm(Isomorphisme de Dolbeault)}
Soit $X$ une vari\'et\'e analytique
complexe, ${\cal F}$ un ${\cal O}_X$-module localement libre sur $X$, 
et $\big(L^{0,q}_{(2), {\rm loc}}(\widetilde
X,\widetilde{\cal F}),\overline\partial\big)$ le complexe des
$(0,q)$-formes sur $\widetilde X$, \`a valeurs dans $\widetilde{\cal
F}$ et localement $L^2$ sur $X$ ainsi que leur $\overline\partial$.  
Alors la cohomologie de ce complexe s'identifie \`a 
$H^\ast_{(2)}(\widetilde X,\widetilde{\cal F})$.
\bigskip

Comme dans la situation classique, l'isomorphisme de Dolbeault fournit 
un moyen commode pour prouver des th\'eor\`emes d'annulations.

\proclaim {\bf 3.6.} Th\'eor\`eme. Soit ${\cal F}$
un faisceau coh\'erent sur un espace complexe~$X$.
\smallskip
\item{\rm a)} Si $U$ est un ouvert de Stein, alors
$H^q_{(2)}(\widetilde U,\widetilde{\cal F})=
H^q\big(U,p_{\ast(2)}\widetilde{\cal F}\big)= 0$
pour  $q > 0$.
\smallskip
\item{\rm b)} Soit $W\subset U$ une paire de Runge
d'ouverts de Stein, alors le morphisme de
restriction
$H^0_{(2)}(\widetilde U,\widetilde{\cal F})\to
H^0_{(2)}(\widetilde W,\widetilde{\cal F})$
est d'image dense pour la topologie de la convergence $L^2$  au dessus
des compacts de~$W$.
\bigskip

\noi {\bf D\'emonstration.} Rappelons qu'une paire d'ouverts de Stein
$W\subset U$ est dite de Runge si l'enveloppe holomorphe convexe de
toute partie compacte de $W$ relativement \`a l'alg\`ebre des
fonctions holomorphes sur $U$ est compacte dans $W$. On sait alors que
l'image du morphisme de restriction ${\cal O}(U)\to{\cal O}(W)$ est
dense, et que pour tout compact $K$ de $W$ il existe une fonction
d'exhaustion strictement plurisousharmonique $\varphi_K$ sur $X$ telle
que l'ensemble de niveau $X_c=\{\varphi_K<c\}$ v\'erifie $K\subset
X_c\compact W$. Le th\'eor\`eme 3.5 se prouve en 3 \'etapes.
\medskip

\noi {\it \'Etape 1. $X$ est lisse et ${\cal F}$ est localement libre.} 
\smallskip

\noi Dans ce cas, il suffit d'utiliser l'isomorphisme de Dolbeault
3.4, et d'appliquer le th\'eor\`eme de
H\"ormander-Andreotti-Vesen\-tini ([H], [AV]) avec des poids 
plurisousharmoniques
$\varphi$ sur $U$ \`a croissance arbitrairement grande lorsqu'on
s'approche du bord de~$U$, de mani\`ere \`a faire converger les
normes~$L^2$. L'assertion sur les paires de Runge se d\'emontre comme
les th\'eor\`emes 4.3.2 et 5.2.10 de H\"ormander [H], en utilisant 
des poids de la forme $e^{-N\varphi_K}$, $N\gg 0$, pour assurer la
convergence uniforme des approximations au voisinage de~$K$. Bien
entendu, cette \'etape permet aussi de couvrir le cas o\`u $U$ est un
ouvert de Stein dans un espace $X$ quelconque, il suffit de plonger
$V$ dans un espace de Stein ambiant lisse et de prolonger le faisceau
${\cal F}$ par $0$ en dehors de $X$.  \medskip

\noi {\it \'Etape 2. L'ouvert $U$ est contenu dans un ouvert de Stein
simple $V$ sur lequel ${\cal F}$ admet une r\'esolution libre.}
\smallskip

\noi Soit
$$
0\to{\cal L}_n\to{\cal L}_{n-1}\to{\cal L}_n\to
\ldots\to{\cal L}_0\to{\cal F}
$$ une r\'esolution libre de ${\cal F}$. On raisonne par r\'ecurrence
sur la longueur $n$ de la r\'esolution.  Si $n=0$, alors ${\cal F}$
est libre et on applique l'\'etape~1. En g\'en\'eral, soit ${\cal G}$
le noyau de ${\cal L}_0\to{\cal F}$. Alors ${\cal G}$ admet une
r\'esolution libre de longueur $n-1$, et par hypoth\`ese de
r\'ecurrence on a $H^q_{(2)}(\widetilde U,\widetilde{\cal G})=0$ pour
$q>0$. La suite exacte
$$
0\to{\cal G}\to{\cal L}_0\to{\cal F}\to 0
$$
fournit une suite exacte longue de cohomologie
$$
0=H^q_{(2)}(\widetilde U,\widetilde{\cal L}_0)
\to H^q_{(2)}(\widetilde U,\widetilde{\cal F})
\to H^{q+1}_{(2)}(\widetilde U,\widetilde{\cal G})
=0,\qquad q>0,$$
ce qui conclut la r\'ecurrence. Le fait que le morphisme 
$H^0_{(2)}(\widetilde W,\widetilde{\cal L}_0)\to
H^0_{(2)}(\widetilde W,\widetilde{\cal F})$
soit surjectif ram\`ene l'assertion sur les paires de Runge au
cas d'un faisceau localement libre.
\medskip

\noi {\it \'Etape 3. Cas g\'en\'eral.}
\smallskip

\noi On utilise la classique ``m\'ethode des bosses'' d'Andreotti-Grauert.
Pour cela, on choisit un recouvrement localement fini 
$U=\bigcup_{j\in{\Bbb N}}U_j$ assez fin de $U$, par des ouverts 
$U_j$ ayant les propri\'et\'es suivantes:

\medskip
\noi$\bullet$\quad $U_j$ est un ouvert de Stein relativement compact
dans $U$, et $(U_j,U)$ est une paire de Runge~;

\smallskip
\noi$\bullet$\quad $V_j=U_0\cup U_1\cup\ldots\cup U_j$ est 
un ouvert de Stein et $(V_j,U)$
est une paire de Runge.
\smallskip

\noi On choisit le recouvrement $(U_j)$ assez fin pour que chaque
$U_j$ soit contenu dans un ouvert de Stein simple sur lequel ${\cal F}$
admet une r\'esolution libre. On d\'emontre maintenant par r\'ecurrence 
sur $j$ que
\smallskip
\item{$(\hbox{a}_j)$}
$H^q_{(2)}(\widetilde V_j,\widetilde{\cal F})=0$ pour tout $q>0$,
\smallskip

\item{$(\hbox{b}_j)$} Si $W\subset V_j$ est une paire
de Runge, alors la restriction $H^0_{(2)}(\widetilde V_j,
\widetilde{\cal F})\to H^0_{(2)}(\widetilde W,\widetilde{\cal F})$
est d'image dense.
\smallskip

\noi
Pour $j=0$ on a $V_0=U_0$ et $(\hbox{a}_0)$, $(\hbox{b}_0)$ r\'esultent
de l'\'etape 2. Pour passer de l'\'etape $j$ \`a l'\'etape $j+1$, on
utilise la suite exacte
$$
\eqalign{\cdots
H^{q-1}_{(2)}(\widetilde V_j, \widetilde{\cal F})&{}\oplus
H^{q-1}_{(2)}(\widetilde U_{j+1}, \widetilde{\cal F})\to
H^{q-1}_{(2)}(\widetilde V_j\cap\widetilde U_{j+1}, \widetilde{\cal F})\to\cr
H^q_{(2)}(\widetilde V_{j+1}, \widetilde{\cal F})
\to H^q_{(2)}(\widetilde V_j, \widetilde{\cal F})
&{}\oplus H^q_{(2)}(\widetilde U_{j+1}, \widetilde{\cal F})~\,\to
H^q_{(2)}(\widetilde V_j\cap\widetilde U_{j+1},\widetilde{\cal F})~\,
\to\cdots\cr}
$$
qui r\'esulte de l'application de la suite exacte de
Mayer-Vietoris au faisceau $(p_2)_\ast\widetilde{\cal F}$. Pour $q\geq 2$,
l'\'etape 2 et l'hypoth\`ese de r\'ecurrence entra\^{\i}nent
$$H^{q-1}_{(2)}(\widetilde V_j\cap\widetilde U_{j+1}, \widetilde{\cal F})
= H^q_{(2)}(\widetilde U_{j+1}, \widetilde{\cal F})=0,\quad\hbox{resp.}\quad
H^q_{(2)}(\widetilde V_j, \widetilde{\cal F})=0,$$ 
d'o\`u $H^q_{(2)}(\widetilde V_{j+1}, \widetilde{\cal F})=0$. Si $q=1$, on
utilise de plus le fait que la restriction
$H^0_{(2)}(\widetilde V_j, \widetilde{\cal F})
\to H^0_{(2)}(\widetilde V_j\cap \widetilde U_{j+1},
\widetilde{\cal F})$ est d'image dense pour voir que le morphisme
continu
$$
H^0_{(2)}(\widetilde V_j\cap
\widetilde U_{j+1}, \widetilde{\cal F})\to
H^1_{(2)}(\widetilde V_{j+1}, \widetilde{\cal F})$$
est n\'ecessairement nul. Ceci implique alors 
$H^1_{(2)}(\widetilde V_{j+1}, \widetilde{\cal F})=0$
et l'assertion $(\hbox{a}_{j+1})$ est d\'emontr\'ee.

L'assertion $(\hbox{b}_{j+1})$, quant \`a elle, s'obtient comme suit.
Soit $W\subset V_{j+1}$ une paire de Runge. Alors $W\cap V_j\subset
V_j$ et $W\cap U_{j+1}\subset U_{j+1}$ sont des paires de Runge pour
lesquelles on peut appliquer l'hypoth\`ese de r\'ecurrence
$(\hbox{b}_j)$, resp.\ l'\'etape 2. Si $s$ est une section de
$H^0_{(2)}(\widetilde W, \widetilde{\cal F})$, on peut approximer $s$
en topologie $L^2$ au dessus de tout compact de $W\cap V_j$, resp.\ de
$W\cap U_{j+1}$, par des sections $s_j\in H^0_{(2)}(\widetilde V_j,
\widetilde{\cal F})$, resp.\ $t_{j+1}\in H^0_{(2)}(\widetilde U_{j+1},
\widetilde{\cal F})$. La diff\'erence $s_j-t_{j+1}$ d\'efinit un
$1$-cocycle de \v Cech sur $V_{j+1}$ relativement au recouvrement
$(V_j,U_{j+1})$.  Comme $H^1_{(2)}(\widetilde V_{j+1},\widetilde{\cal F})
=0$, on a un morphisme surjectif d'espaces de Fr\'echet
$$
H^0_{(2)}(\widetilde V_j, \widetilde{\cal F})\oplus
H^0_{(2)}(\widetilde U_{j+1}, \widetilde{\cal F})\to
H^0_{(2)}(\widetilde V_{j+1}, \widetilde{\cal F}).
$$
Or, si $s_j$ et $t_{j+1}$ sont des approximations suffisamment
bonnes de $s$, la diff\'erence $s_j-t_{j+1}$ peut \^etre choisie
arbitrairement petite dans la topologie de l'espace de Fr\'echet but.
D'apr\`es le th\'eor\`eme de l'application ouverte, on peut trouver
des sections $\sigma_j\in H^0_{(2)}(\widetilde V_j,\widetilde{\cal F})$
et \hbox{$\tau_{j+1}\in H^0_{(2)}(\widetilde U_{j+1},\widetilde{\cal F})$}
arbitrairement petites telles que
$\sigma_j-\tau_{j+1}=s_j-t_{j+1}$ sur $\widetilde V_j\cap\widetilde
U_{j+1}$. Alors $s_j-\sigma_j$ et $t_{j+1}-\tau_{j+1}$ se recollent en
une section sur $V_{j+1}=V_j\cup U_{j+1}$ qui approxime $s$ d'aussi
pr\`es qu'on veut sur~$W$.  Un raisonnement standard de passage \`a la
limite \`a la Mittag-Leffler permet d'atteindre la nullit\'e de la
cohomologie sur $U=\bigcup V_j$ et le th\'eor\`eme de Runge pour la
paire $W\subset U$ \`a partir du th\'eor\`eme de Runge sur les paires
$W\cap V_j\subset V_{j+1}$, $V_{j+1}\subset V_{j+2}$, $\ldots\,$, etc.
\bigskip

\proclaim {\bf 3.7.} Corollaire. Soit ${\cal U} :=
(U_\lambda)_{\lambda\in \Lambda}$ un recouvrement ouvert localement
fini de $X$ par des ouverts de Stein $U_\lambda$. Ce recouvrement est 
alors de Leray pour $p_{\ast(2)}{\cal F}$ et on a
un isomorphisme naturel :
$$
H^\ast_{(2)}(\widetilde {\cal U},\widetilde{\cal F}) 
\buildrel \sim \over {\to}H^\ast_{(2)}(\widetilde X,\widetilde{\cal F}),
$$
o\`u $H^\ast_{(2)}(\widetilde {\cal U},\widetilde{\cal F}) := H^\ast({\cal
U},p_{\ast(2)}\widetilde{\cal
F})$ est la cohomologie de \v Cech de $p_{\ast(2)}\widetilde{\cal F}$
relative au
recouvrement ${\cal U}$ de $X$.
\bigskip

\noi {\bf 3.8.} Explicitons l'assertion de 3.7 : soit $N_q({\cal U})$
le $q$-nerf du recouvrement ${\cal U}$, constitu\'e des
intersections non vides de $(q+1)$ des \'el\'ements de ${\cal U}$, et
soit $C^q_{(2), {\rm loc}}(\widetilde {\cal U},\widetilde{\cal F})$ le 
groupe des $q$-cocha\^{\i}nes \`a valeurs dans $\widetilde{\cal F}$,
d\'efinies sur les ouverts $\widetilde{\cal U}_{(q)}=
p^{-1}({\cal U}_{(q)})$ et localement $L^2$ au dessus de ${\cal U}_{(q)}
\in N_q({\cal U})$, avec les applications de cobord 
$$
\delta_q : C^q_{(2), {\rm loc}}(\widetilde 
{\cal U},\widetilde{\cal F}) \to C^{q+1}_{(2), {\rm loc}}(\widetilde 
{\cal U},\widetilde{\cal F})
$$ 
usuelles. Alors $H^\ast_{(2)}(\widetilde {\cal U},\widetilde{\cal F})$ 
est la cohomologie du complexe ainsi d\'efini. Les espaces
$C^q_{(2), {\rm loc}}(\widetilde {\cal U},\widetilde{\cal F})$ 
sont naturellement munis de la topologie de la convergence en norme
$L^2$ au dessus des compacts $p$-simples contenus dans les  ${\cal U}_{(q)}$
(on ne prend bien entendu en compte simultan\'ement qu'un nombre fini de
ces intersections), ce qui en fait des espaces de Fr\'echet. On munit
$H^\ast_{(2)}(\widetilde {\cal U},\widetilde{\cal F})$ de la topologie
quotient correspondante (qui n'est pas n\'ecessairement s\'epar\'ee).

\bigskip

\noi {\bf 3.9.} D\'esignons par $B^q_{(2)}(\widetilde {\cal U},
\widetilde{\cal F})$ et $Z^q_{(2)}(\widetilde {\cal U},\widetilde{\cal
F})$ respectivement l'image de $\delta_{q-1}$ et le noyau de $\delta_q$
$\big($avec $B^0_{(2)}(\widetilde {\cal U},\widetilde{\cal F})=0\big)$. 
On a alors une suite exacte
$$
0 \to \underline H^q_{(2)}(\widetilde {\cal U},\widetilde{\cal F})
\to H^q_{(2)}(\widetilde {\cal U},\widetilde{\cal F}) \to
\overline H^q_{(2)}(\widetilde {\cal U},\widetilde{\cal F})
\to 0,
$$
o\`u $\overline B^q_{(2)}(\widetilde {\cal U},\widetilde{\cal F})$ est 
l'adh\'erence dans $Z^q_{(2)}(\widetilde {\cal U},
\widetilde{\cal F})$ de $B^q_{(2)}(\widetilde {\cal U},\widetilde{\cal F})$
et
$$
\underline H^q_{(2)}(\widetilde {\cal U},\widetilde{\cal F})=
\overline B^q_{(2)}(\widetilde {\cal U},\widetilde{\cal F})/
B^q_{(2)}(\widetilde {\cal U},\widetilde{\cal F}),\qquad
\overline H^q_{(2)}(\widetilde {\cal U},\widetilde{\cal F})
=Z^q_{(2)}(\widetilde {\cal U},\widetilde{\cal F})
/\overline B^q_{(2)}(\widetilde {\cal U},\widetilde{\cal F}).
$$
L'espace $\overline H^q_{(2)}(\widetilde {\cal U},\widetilde{\cal F})$
est par d\'efinition un espace de Fr\'echet, mais
$\underline H^q_{(2)}(\widetilde {\cal U},\widetilde{\cal F})$
est muni de la topologie grossi\`ere et, s'il est non nul, la cohomologie 
$L^2$ n'est pas s\'epar\'ee. On va voir, cependant, que la topologie 
de $H^q_{(2)}(\widetilde {\cal U},\widetilde{\cal F})$ est essentiellement 
ind\'ependante du choix du recouvrement. Soient en effet ${\cal U}'$, 
${\cal U}$ des recouvrements de Stein de $X$. On suppose ${\cal U}'$ plus 
fin que  ${\cal U}$ et muni d'une application de raffinement $\rho$ vers 
${\cal U}$. Alors, dans le diagramme commutatif associ\'e $(q \ge 0)$,
$$\matrix{
0 &\hfl{}{} &\underline H^q_{(2)}(\widetilde {\cal U},\widetilde{\cal F}) 
&\hfl{}{}
&H^q_{(2)}(\widetilde {\cal
U},\widetilde{\cal F}) &\hfl{}{} &\overline H^q_{(2)}(\widetilde {\cal
U},\widetilde{\cal F}) &\hfl{}{} &0
\cr
\noalign{\vskip 2mm}
&&\vfl{}{\underline\rho} &&\vfl{}{\rho} &&\vfl{}{\overline\rho} \cr
\noalign{\vskip 2mm}
0 &\hfl{}{} &\underline H^q_{(2)}(\widetilde {\cal U}',\widetilde{\cal F}) 
&\hfl{}{} &H^q_{(2)}(\widetilde {\cal
U}',\widetilde{\cal F}) &\hfl{}{} &\overline H^q_{(2)}(\widetilde {\cal
U}',\widetilde{\cal F}) &\hfl{}{} &0
\cr}$$
les applications verticales $\rho$, $\underline \rho$, $\overline\rho$
sont des isomorphismes topologiques (d'espaces de Fr\'echet en ce qui
concerne $\overline\rho$). En effet, si (pour simplifier) on d\'esigne
par $(C,\delta)$ et $(C',\delta')$ les complexes de Fr\'echet
impliqu\'es, l'isomorphisme de Leray implique que $\rho$ est un
isomorphisme alg\'ebrique, et il est clair par ailleurs que
l'application de restriction $\rho:C\to C'$ est continue.  On a alors
une application surjective $\delta'\oplus r:C'\oplus Z\to Z'$ entre
espaces de Fr\'echet. Le th\'eor\`eme de l'application ouverte montre
que cette application est ouverte, et il en est donc de m\^eme pour
l'application induite $\rho:H=Z/\delta C\to H'=Z'/\delta'C'$. Ceci
montre d\'ej\`a que $\rho$ est un isomorphisme
topologique. L'assertion pour $\underline\rho$ et $\overline\rho$ s'en
d\'eduit imm\'ediatement, puisque $\underline H$ est la partie
grossi\`ere et $\overline H$ le quotient s\'epar\'e de la cohomologie.
Ceci nous m\`ene \`a la d\'efinition suivante.
\bigskip

\proclaim {\bf 3.10.} D\'efinition. L'espace de cohomologie 
$H^q_{(2)}(\widetilde X,\widetilde{\cal F})$ $(q\ge 0)$ est muni
d'une topologie naturelle pour laquelle, si 
$\underline H^q_{(2)}(\widetilde X,\widetilde{\cal F})$ est l'adh\'erence de
$0$, alors 
$$
\overline H^q_{(2)}(\widetilde X,\widetilde{\cal F})=
H^q_{(2)}(\widetilde X,\widetilde{\cal F})/\underline H^q_{(2)}
(\widetilde X,\widetilde{\cal F})
$$ 
est un espace de Fr\'echet. On appellera
$\overline H^q_{(2)}(\widetilde X,\widetilde{\cal F})$ $($resp.\ $\underline
H^q_{(2)}(\widetilde X,\widetilde{\cal F}))$ la {\it cohomologie $L^2$ 
s\'e\-pa\-r\'ee} de $X$ \`a valeurs dans $\widetilde{\cal F}$
$($resp.\ le \hbox{\it noyau} de la cohomologie $L^2$ de $X$ \`a valeurs 
dans~$\widetilde{\cal F}\,)$.
\bigskip

\noi {\bf 3.11. Remarque.} Par construction, les objects $\underline
H^q_{(2)}(\widetilde X,\widetilde{\cal F})$
et $\overline H^q_{(2)}(\widetilde X,\widetilde{\cal F})$ sont fonctoriels 
en ${\cal F}$ et $X$.
\bigskip

\noi {\bf 3.12. Remarque.} Il r\'esulte imm\'ediatement de 2.9 que si
${\cal F}$ est support\'e par la sous-vari\'et\'e $Y$ de $X$, alors :
$H^\ast_{(2)}(\widetilde X,\widetilde{\cal F}) \mapsto
H^\ast_{(2)}(\widetilde Y,\widetilde{\cal F})$ est un isomorphisme
topologique.
\bigskip

\noi {\bf 3.13. Remarque.} Si $X$ est un espace compact, on
peut choisir un recouvrement ouvert fini ${\cal U}=(U_j)$ fini par des ouverts
de Stein $p$-simples, puis r\'etr\'ecir un peu chacun des ouverts $U_j$
en des ouverts $U''_j\compact U'_j\compact U_j$ tels que
$U'_j$ et $U''_j$ soient de Runge dans $U_j$. Les morphismes de restriction
donnent lieu \`a des fl\`eches
$$
C^q_{(2), {\rm loc}}(\widetilde {\cal U},\widetilde{\cal F})
\to
C^q_{(2)}(\widetilde {\cal U'},\widetilde{\cal F})
\to
C^q_{(2), {\rm loc}}(\widetilde {\cal U''},\widetilde{\cal F})
$$
o\`u les termes extr\^emes sont des espaces de Fr\'echet et le terme central 
un espace de Hilbert (on prend sur ce terme la topologie $L^2$ globale
sur ${\cal U'}$). En cohomologie, on un isomorphisme entre les termes
extr\^emes, ce qui prouve que la cohomologie du terme central se surjecte
sur cette cohomologie. La cohomologie s\'epar\'ee $\overline H^q_{(2)}
(\widetilde X,\widetilde{\cal F})$ poss\`ede donc alors une topologie d'espace
de Hilbert.

\bigskip

\noi {\bf 3.14. Remarque.} Supposons maintenant que $X$ soit une
vari\'et\'e compacte lisse et que ${\cal F}$ soit un faisceau
analytique localement libre sur $X$.  Les arguments de [G]
s'appliquent encore dans ce contexte et montrent que si
$Z^{0,q}_{(2)}(\widetilde X,\widetilde{\cal F}) \subset
C^{0,q}_{(2)}(\widetilde X,\widetilde{\cal F})$ est le noyau du
$\overline\partial$, alors : $\big(Z^{0,q}_{(2)}(\widetilde
X,\widetilde{\cal F})/(\overline{\hbox{Im} \,\overline\partial})\big)$
s'identifie \`a l'espace de Hilbert (par ellipticit\'e de $\overline
\partial$) des formes $\Delta_{\overline\partial}$-harmoniques de type
$(0,q)$ et $L^2$ sur $\widetilde X$ \`a valeurs dans $\widetilde{\cal
F}$. On voit, de plus, que ${\cal H}^{0,q}_{(2)}(\widetilde
X,\widetilde{\cal F})$ s'identifie canoniquement \`a la cohomologie
r\'eduite $\overline H^q_{(2)}(\widetilde X,\widetilde{\cal F})$
d\'efinie en 3.10.
\bigskip

\proclaim {\bf 3.15.} Corollaire \hbox{\rm (dualit\'e de Serre)}. Soit
$X$ une vari\'et\'e complexe compacte lisse, et ${\cal F}$ un faisceau
analytique coh\'erent localement libre sur $X$. Soit $p : \widetilde X
\to X$ un rev\^etement \'etale. Il existe une isom\'etrie
antilin\'eaire $\sigma : \overline H^q_{(2)}(\widetilde
X,\Omega^r_{\widetilde X} \otimes \widetilde{\cal F}) \to \overline
H^{n-q}_{(2)}(\widetilde X,\Omega^{n-r}_{\widetilde X} \otimes
\widetilde{\cal F}^\ast)$ si $n$ est la dimension $($pure$)$ de $X$.
\bigskip

\noi {\bf 3.16. Cas galoisien.} Dans le cas particulier o\`u le
rev\^etement $p : \widetilde X \to X$ de la vari\'et\'e complexe
compacte $X$ est galoisien, de groupe $\Gamma$, on a une op\'eration
naturelle du groupe $\Gamma$ sur tous les objets d\'efinis
pr\'ec\'edemment : $\widetilde{\cal F}$, $p_{\ast(2)}\widetilde{\cal
F}$, $C^q_{(2),{\rm loc}}(\widetilde {\cal U},\widetilde{\cal F})$,
$H^q_{(2)}(\widetilde {\cal U},\widetilde{\cal F})$, $\ldots\;$, et ce
pour tout $q\ge 0$ et tout recouvrement de Stein $p$-simple ${\cal
U}$. On a donc aussi une action de $\Gamma$ sur les espaces de
cohomologie $H^q_{(2)}(\widetilde X,\widetilde{\cal F})$, $\underline
H^q_{(2)}(\widetilde X,\widetilde{\cal F})$, $\overline
H^q_{(2)}(\widetilde X,\widetilde{\cal F})$.  Dans le cas o\`u $X$ est
compacte, cette op\'eration induit une action unitaire sur l'espace de
Hilbert $C^q_{(2)}(\widetilde {\cal U},\widetilde{\cal F})$.
\bigskip

\proclaim {\bf 3.17.} Proposition. Soit $p : \widetilde X \to X$ un
rev\^etement galoisien de groupe $\Gamma$ de la vari\'et\'e complexe
compacte $X$ ; soit ${\cal F}$ un faisceau coh\'erent sur $X$, et
${\cal U}$ un recouvrement de Stein ouvert, fini et $p$-simple de
$X$. Cette action d\'efinit une action de $\Gamma$ sur
$H^\ast_{(2)}(\widetilde X,\widetilde{\cal F})$ qui pr\'eserve chacune
des semi-normes pr\'ehilbertiennes (\'equivalentes entre elles) dont 
cet espace peut \^etre muni. En particulier, $\Gamma$ agit sur $\underline
H^\ast_{(2)}(\widetilde X,\widetilde{\cal F})$ et $\overline
H^\ast_{(2)}(\widetilde X,\widetilde{\cal F})$, de mani\`ere unitaire
sur ce dernier espace $($qui est de Hilbert$)$.
\vskip 0.75cm

\noi {\bigbf \S\ 4. Th\'eor\`emes d'annulation.}
\bigskip

\noi
{\bf 4.0.} Les nombreux r\'esultats d'annulation accessibles par les
techniques $L^2$ usuelles vont en g\'en\'eral se transcrire mot pour mot
pour donner des versions s'appliquant en cohomologie~$L^2$. Nous
pr\'esentons ici quelques \'enonc\'es parmi les plus fondamentaux.
\bigskip

\proclaim {\bf 4.1.} Th\'eor\`eme de Kodaira-Serre $L^2$. Soit 
$X$ une vari\'et\'e projective
lisse, et ${\cal F}$ un faisceau analytique coh\'erent sur~$X$. 
Soit ${\cal L}$ un fibr\'e en droites ample sur $X$ et $p :
\widetilde X \to X$ un rev\^etement \'etale de $X$. Il existe $m_0 =
m_0({\cal L},{\cal F})$, ind\'ependant du rev\^etement $p$, tel que
$H^q_{(2)}\big(X,\widetilde{\cal F}(m)\big) = 0$ si $q > 0$ et $m
\ge m_0$ $($on pose ici comme d'habitude ${\cal F}(m) := {\cal F} 
\otimes {\cal L}^m)$.
\bigskip

\noi {\bf D\'emonstration.} ${\cal F}$ admet une r\'esolution
localement libre de longueur $0 \le r \le n = \dim_{\Bbb C}(X)$. Si
$r=0$, ${\cal F}$ est localement libre, et le r\'esultat est
cons\'equence directe de 3.5 ci-dessus et de l'existence de solutions
$L^2$ \`a l'\'equation : $\overline \partial \tilde v = \tilde w$,
avec $\overline\partial\tilde w=0$ et $\tilde w$ section $L^2$ de
$\widetilde{\cal F}^{0,1}$ ([D], th\'eor\`eme 5.1).

\noi Sinon, on proc\`ede \`a nouveau par r\'ecurrence sur $r$,
supposant le r\'esultat vrai pour $r-1 \ge 0$. Dans ce cas,
l'assertion r\'esulte imm\'ediatement de la suite exacte longue de
cohomologie $L^2$ (th\'eor\`eme 3.2) associ\'ee \`a la suite exacte de
faisceaux : $0 \to {\cal G} \to {\cal H} \to {\cal F} \to 0$, o\`u
${\cal H}$ est localement libre et o\`u ${\cal G}$ admet une
r\'esolution localement libre de longueur $(r-1)$.\cqfd
\bigskip

\noi {\bf 4.2. Exemple} (Cet exemple a partiellement motiv\'e la
construction pr\'esent\'ee ici). Soit $X$ une vari\'et\'e projective,
${\cal F}$ et ${\cal L}$ des faisceaux analytiques coh\'erents sur
$X$, avec ${\cal L}$ fibr\'e en droites ample. Soit $p : \widetilde X 
\to X$ un rev\^etement \'etale de $X$ et $Y$ un sous-sch\'ema (non
n\'ecessairement r\'eduit de $X$). Il existe
$m_0 := m_0({\cal F},{\cal L},X,Y)$ tel que, pour $m \ge m_0$
le morphisme de restriction naturel
$$
H^0_{(2)}(\widetilde X,\widetilde{\cal F} \otimes \widetilde{\cal L}^m)
\longrightarrow
H^0_{(2)}(\widetilde Y,\widetilde{\cal F}_{|Y} \otimes \widetilde
{\cal L}_{|Y}^m)
$$
soit surjectif. On utilise en effet la suite-exacte
$$
0\to{\cal I}_Y{\cal F}\to{\cal F}\to{\cal F}_{|Y}\to 0
$$
o\`u ${\cal F}_{|Y}={\cal F}/{\cal I}_Y{\cal F}$ est la restriction de
$\cal F$ au sous-sch\'ema $Y$. Le th\'eor\`eme de Kodaira-Serre implique
l'annulation du groupe
$$
H^1_{(2)}(\widetilde X,\widetilde{\cal I}_{\widetilde Y}
\widetilde{\cal F}\otimes \widetilde{\cal L}^m)
$$
pour $m\ge m_0$ assez grand, d'o\`u le r\'esultat. Ceci s'applique entre
autres au cas o\`u $Y$ est le sch\'ema ponctuel associ\'e \`a l'anneau 
${\cal O}_X/{\cal I}_a^{k+1}$ des jets d'ordre $k$ de fonctions en 
un point $a$ de $X$. On voit alors que les $k$-jets
de $p_{*(2)}\big(\widetilde{\cal F} \otimes \widetilde{\cal L}^m\big)$ sont
engendr\'es pour $m\ge m_0$ assez grand par les sections 
globales $L^2$ du faisceau $\widetilde{\cal F} \otimes \widetilde{\cal
L}^m$ sur $\widetilde X$ (il faut voir que $m_0$ peut \^etre choisi
ind\'ependant de $a$, mais c'est imm\'ediat en contr\^olant un tant
soit peu les r\'esolutions libres globales des anneaux ${\cal O}_X/
{\cal I}_a^{k+1}$).\cqfd 
\medskip

Les r\'esultats suivants sont des transcriptions imm\'ediates des r\'esultats
$L^2$ classiques pour les vari\'et\'es k\"ahl\'eriennes compl\`etes
([AV], [D]), et nous les \'enon\c{c}ons donc sans commentaires
(le corollaire 4.5 \'etant par exemple d\'ej\`a mentionn\'e dans [K], 11.4).
\bigskip

\proclaim {\bf 4.3.} Th\'eor\`eme de Akizuki-Kodaira-Nakano $L^2$. Soit 
$X$ une vari\'et\'e projective lisse de dimension complexe $n$,
$p : \widetilde X \to X$ un rev\^etement \'etale de $X$ et 
${\cal L}$ un fibr\'e en droites ample sur $X$. Alors on a 
$$
H^q_{(2)}\big(\widetilde X,\Omega^r_{\widetilde X}\otimes
\widetilde{\cal L}\big) = 0
\qquad 
\hbox{si $q+r\ge n+1$.}
$$
\bigskip

\proclaim {\bf 4.4.} Th\'eor\`eme de Nadel $L^2$. \kern-1pt Soit 
$X$ une vari\'et\'e compacte $($projective ou de Moishezon$)$, lisse,
$p : \widetilde X \to X$ un rev\^etement \'etale de $X$ et 
${\cal L}$ un fibr\'e en droites sur $X$. On suppose que
${\cal L}$ poss\`ede une m\'etrique hermitienne singuli\`ere $h$
dont la $(1,1)$-forme de courbure $\Theta_h({\cal L})$ est positive
au sens des courants, minor\'ee par une $(1,1)$ forme de classe
$C^\infty$ d\'efinie positive. Alors on a
$$
H^q_{(2)}\big(\widetilde X,K_{\widetilde X}\otimes\widetilde{\cal L}
\otimes\widetilde{\cal I}(h)\big) = 0
\qquad 
\hbox{pour tout $q\ge 1$,}
$$
o\`u ${\cal I}(h)\subset{\cal O}_X$ d\'esigne l'id\'eal multiplicateur 
des germes de fonctions holomorphes $f$ telles que $\int|f|^2e^{-\varphi}
<+\infty$ $(e^{-\varphi}$ d\'esignant le poids qui repr\'esente localement la
m\'etrique $h)$.
\bigskip

\proclaim {\bf 4.5.} Corollaire (Th\'eor\`eme de Kawamata-Viehweg $L^2$).
Si $X$ est une vari\'et\'e de Moishezon lisse et
$p : \widetilde X \to X$ un rev\^etement \'etale de $X$. On
suppose donn\'e un fibr\'e en droites ${\cal L}$ num\'eriquement \'equivalent 
\`a la somme d'un $\bQ$-diviseur $D$ nef $($num\'eriquement effectif$)$,
et d'un $\bQ$-diviseur effectif $E$.
\smallskip
\item{\rm(i)}Si $D$ est gros, alors
$H^q_{(2)}\big(\widetilde X,K_{\widetilde X}\otimes\widetilde{\cal L}
\otimes\widetilde{\cal I}(E)\big) = 0$
pour tout $q\ge 1$.
\smallskip
\item{\rm(ii)} Si $D$ est de dimension num\'erique $\nu\le n=\dim X$,
l'annulation a lieu pour $q>n-\nu$.
\smallskip
\item{\rm(iii)} Si $D$ est de dimension num\'erique $\nu\le n$
et si l'id\'eal ${\cal I}(E)$ est un faisceau inversible $($i.e.\ l'id\'eal
d'un diviseur effectif$)$, la cohomologie 
s\'epar\'ee
$\overline H^q_{(2)}\big(\widetilde X,\widetilde{\cal L}^{-1}
\otimes(\widetilde{\cal I}(E))^{-1}\big)$ est nulle pour tout $q<\nu$.
\bigskip

\noi {\bf D\'emonstration.} 
Rappelons que ${\cal I}(E)$ est le faisceau associ\'e \`a $\varphi={1\over k}
\log|g|$, o\`u $g$ est un g\'en\'erateur de ${\cal O}(-kE)$.
La preuve de 4.5 consiste en une r\'eduction \`a 4.4, \`a peu pr\`es
identique \`a celle effectu\'ee dans le cas classique. 
\smallskip

\noi
(i) Si $D$ est ample, le r\'esultat r\'esulte directement de 4.4, en 
munissant ${\cal O}(D)$ d'une m\'etrique lisse \`a courbure positive et
${\cal O}(E)$ de la m\'etrique associ\'ee au poids $e^{-\varphi}$ (dont la
courbure est le courant d'int\'egration $[E]$). En g\'en\'eral, si $D$ est 
seulement nef et gros, on peut \'ecrire
$D=D'+F$ avec $D'$ ample et $F$ un $\bQ$-diviseur effectif aussi petit
que l'on veut. On peut en particulier supposer que 
${\cal I}(E+F)={\cal I}(E)$. 
\smallskip

\noi
(ii) On se ram\`ene au cas o\`u la dimension num\'erique est maximale par un
argument classique de sections hyperplanes et un raisonnement par r\'ecurrence
sur la dimension. De fa\c{c}on pr\'ecise, on choisit un diviseur lisse
$Y$ tr\`es ample dans $X$ et on consid\`ere la suite exacte courte
$$
0\to K_X\otimes{\cal L}\to K_X\otimes{\cal O}(Y)\otimes{\cal L}\to
K_Y\otimes{\cal L}_{|Y}\to 0.
$$
L'annulation de la cohomologie du terme central est obtenue par 
Kodaira-Serre en prenant $Y$ assez grand, tandis que l'annulation de
la cohomologie en degr\'e $q-1$ du terme de droite r\'esulte de 
l'hypoth\`ese de r\'ecurrence.
\smallskip

\noi
(iii) C'est un cas particulier de (ii), si on utilise la dualit\'e de 
Serre. Il serait int\'eressant de savoir si la cohomologie non
s\'epar\'ee $H^q_{(2)}\big(\widetilde X,\widetilde{\cal L}^{-1}
\otimes(\widetilde{\cal I}(E))^{-1}\big)$
est nulle elle aussi. La difficult\'e est que c'est une cohomologie
``duale'' d'une cohomologie $L^2$, qui ne s'obtient pas directement
par application d'estimations $L^2$ globales.\cqfd
\vskip0.75cm

\noi {\bigbf \S\ 5. Th\'eor\`eme de finitude et th\'eor\`eme de l'indice.}
\bigskip
\noi {\bf 5.0.}
Notre objectif est ici d'\'etendre au cas de faisceaux analytiques coh\'erents 
quelconques le th\'eor\`eme de l'indice $L^2$ de Atiyah~[A]. La preuve en est
purement formelle \`a partir des r\'esultats des sections pr\'ec\'edentes.
\bigskip

\proclaim {\bf 5.1.} Th\'eor\`eme. Soit $X$ un espace analytique
compact et ${\cal F}$ un faisceau analytique coh\'erent sur $X$. Soit $p :
\widetilde X \to X$ un rev\^etement \'etale galoisien de groupe
$\Gamma$. Pour tout $q \ge 0$, le groupe $H^q_{(2)}(\widetilde
X,\widetilde {\cal F})$ est un $\Gamma$-module $L^2$ de pr\'esentation
finie. En particulier, la $\Gamma$-dimension de $\overline
H^q_{(2)}(\widetilde X,\widetilde{\cal F})$, not\'ee
$h^q_{(2)}(\widetilde X,\widetilde{\cal F})$ est (un nombre r\'eel)
fini. De plus, la caract\'eristique d'Euler $L^2$ 
$$\chi_{(2)}(\widetilde X,\widetilde{\cal F}) := \sum^n_{q=0}(-1)^q
h^q_{(2)}(\widetilde X,\widetilde{\cal F})$$
sur $\widetilde X$ est \'egale \`a la caract\'eristique d'Euler ordinaire:
$$\chi_{(2)}(\widetilde X,\widetilde{\cal F}) = \chi(X,{\cal F}) := 
\sum^n_{q=0}(-1)^q h^q(X,{\cal F}).$$

\noi(Voir l'appendice pour les notions hilbertiennes requises, en particulier 
6.4 et 6.5).
\bigskip

\noi {\bf D\'emonstration.} Si $X$ est lisse et ${\cal F}$ localement 
libre, c'est le th\'eor\`eme de l'indice $L^2$ d'Atiyah ([A]). En 
g\'en\'eral, on raisonne par r\'ecurrence sur la dimension $n=\dim X$, 
en utilisant un d\'evissage de ${\cal F}$ et une r\'esolution des
singularit\'es. Supposons le th\'eor\`eme d\'ej\`a d\'emontr\'e en
dimension $n-1\;$; les r\'esultats sont triviaux en dimension $0$,
car si $X=\{p\}$, on a 
$$\chi(X,{\cal F})=h^0(\{p\},{\cal F})=\dim {\cal F}_p$$ 
tandis que
$H^0_{(2)}(\widetilde X, \widetilde{\cal F})=
\ell^2(\Gamma)\otimes {\cal F}_p$,
d'o\`u $\chi_{(2)} (\widetilde X,\widetilde{\cal F})=\dim{\cal F}_p$.
\medskip

En g\'en\'eral, si $X$ n'est pas r\'eduit, on peut consid\'erer sa
r\'eduction $X_{{\rm red}}$ et la filtration de ${\cal F}$ par les
${\cal N}^p{\cal F}$, o\`u ${\cal N}$ d\'esigne l'id\'eal des
\'el\'ements nilpotents de ${\cal O}_X$.  Le gradu\'e ${\cal N}^p{\cal
F}/{\cal N}^{p+1}{\cal F}$ de cette filtration est constitu\'e de
faisceaux coh\'erents sur ${\cal O}_{X_{{\rm red}}}$.  Par
additivit\'e de la caract\'eristique d'Euler (ordinaire ou $L^2$,
gr\^ace \`a 3.2), on est ramen\'e au cas o\`u $X$ est r\'eduit.
Si $X$ n'est pas lisse, on utilise le th\'eor\`eme de Hironaka pour
trouver une d\'esingu\-larisation $f:Y\to X$. Soit $p':\widetilde Y$
le rev\^etement image r\'eciproque de $p:\widetilde X\to X$ par $f$ et
${\cal G}=f^\ast {\cal F}$, qui est un faisceau coh\'erent sur $Y$. On
a un morphisme naturel injectif ${\cal F}\to f_\ast{\cal G}$, et le
faisceau quotient $f_\ast{\cal G}/{\cal F}$ est \`a support dans le
lieu singulier $X_{{\rm sing}}$. D'apr\`es l'hypoth\`ese de
r\'ecurrence, les r\'esultats sont vrais pour le faisceau quotient
$f_\ast{\cal G}/{\cal F}$, et on est donc ramen\'e \`a traiter le cas
du faisceau image directe $f_\ast{\cal G}$. On utilise alors les suites
spectrales de Leray 
$$
\eqalign{
H^\ast(X,R^\ast f_\ast{\cal G}) &{}\Longrightarrow H^\ast(Y,{\cal G}),\cr
H^\ast_{(2)}(\widetilde X,R^\ast{\widetilde f}_{\ast(2)}{\widetilde{\cal G}}) 
&{}\Longrightarrow H^\ast_{(2)}(\widetilde Y,{\widetilde{\cal G}})\cr}
$$ 
(l'existence de la deuxi\`eme suite spectrale est une cons\'equence
du th\'eor\`eme de Leray et de la Proposition 3.3). Supposons que le
r\'esultat soit d\'ej\`a d\'emontr\'e dans le cas d'une vari\'et\'e
lisse. Alors $\chi(Y,{\cal G})=\chi_{(2)}(\widetilde
Y,{\widetilde{\cal G}})$ puisque $Y$ est lisse, et de m\^eme pour
$q>0$ on a $\chi(X,R^qf_\ast{\cal G})=\chi_{(2)}(\widetilde X,
R^q{\widetilde f}_{\ast(2)}{\widetilde{\cal G}})$ pour $q>0$ par
hypoth\`ese de r\'ecurrence sur la dimension (les faisceaux
$R^qf_\ast{\cal G}$, $q>0$, sont support\'es par $X_{{\rm sing}}$ qui
est de dimension${}\le n-1$).  Comme il y a pr\'eservation de la
caract\'eristique d'Euler dans les diff\'erents niveaux d'une suite
spectrale, il y a \'equivalence \`a prouver le r\'esultat pour la
paire $(X,R^0 f_\ast{\cal G})$ ou pour la paire $(Y,{\cal G})$, ce qui
fait qu'on est ramen\'e au cas o\`u l'espace ambiant $X$ est lisse de
dimension~$n$.

Dans ce cas, soit ${\cal F}_{{\rm tors}}$ la partie de torsion de 
${\cal F}$. Cette partie est \`a support en codimension $n-1$, donc
l'hypoth\`ese de r\'ecurrence s'y applique. La suite exacte 
$$
0\to {\cal F}_{{\rm tors}} \to {\cal F} \to {\cal F}/{\cal F}_{{\rm tors}}
\to 0
$$
et l'additivit\'e de la caract\'eristique d'Euler montre que l'on peut
supposer ${\cal F}$ sans torsion. En appliquant de nouveau le
th\'eor\`eme de Hironaka, il existe une modification analytique propre
$f:Y\to X$, tel que ${\cal G}=f^\ast{\cal F}$ soit localement libre
sur~$Y$. Des arguments identiques \`a ceux qui pr\'ec\`edent
ram\`enent la preuve du cas de la paire $(X,{\cal F})$ au cas de la
paire $(Y,{\cal G})$. On conclut en appliquant cette fois le th\'eor\`eme
de l'indice $L^2$ de Atiyah [A] \`a $(Y,{\cal G})$.\hbox{\cqfd} 
\vskip 1cm

\noi {\bigbf \S\ 6. Appendice : pr\'esentation hilbertienne et
$\Gamma$-dimension.}

\bigskip

\proclaim {\bf 6.1.} D\'efinition. Soit $H$ un espace vectoriel
complexe. Une \hbox{\it pr\'esentation hilbertienne} de $H$ est la
donn\'ee d'une application lin\'eaire continue $\delta : C \to Z$
entre deux espaces de Hilbert et d'un isomorphisme (alg\'ebrique)
$(Z/\delta C) \buildrel \sim \over {\mapsto} H$. On note $(\delta C)$
(resp.  $(\overline{\delta C})$) l'image dans $Z$ de $C$ (resp. son
adh\'erence).\hfill\break 
Associ\'ee \`a une telle pr\'esentation est d\'efinie une suite exacte :
$$
0 \to \underline H =: (\overline{\delta C}/\delta C) \to H =: (Z/\delta
C) \to \overline H:= (Z/\overline{\delta C}) \to 0.
$$
On appelle $\underline H$ (resp.\ $\overline H$) le \hbox{\it noyau} (resp.\
la \hbox{\it r\'eduction}) de $H$ relative \`a cette pr\'esentation 
hilbertienne.\hfill\break
(Ces notions ne d\'ependent que de la classe d'\'equivalence des espaces de
Hilbert).\hfill\break
Une application entre deux pr\'esentations hilbertiennes $\delta : C \to Z$
et $\delta' : C'
\to Z'$ est un diagramme commutatif $s : C \to C'$ et $r : Z \to Z'$ 
d'applications lin\'eaires continues telles que $r\delta=\delta's$. 
Une telle application induit un diagramme commutatif :
$$\matrix{
0 &\hfl{}{} &\underline H &\hfl{}{} &H &\hfl{}{} &\overline H &\hfl{}{} &0 \cr
\noalign{\vskip 2mm}
&&\vfl{}{\underline r} &&\vfl{}{[r]} &&\vfl{}{\overline r} \cr
\noalign{\vskip 2mm}
0 &\hfl{}{} &\underline H' &\hfl{}{} &H' &\hfl{}{} &\overline H' &\hfl{}{}
&0 \cr} \leqno
(6.1')$$
Deux pr\'esentations hilbertiennes sont dites \hbox{\it compatibles} s'il
existe une application entre elles.

\bigskip

\proclaim {\bf 6.2.} Proposition. La situation \'etant celle d\'ecrite en 6.1,
alors :\hfill\break
\noi \hbox{\bf i.} Si $[r]$ est surjective, $\overline r$ est surjective.
Si $[r]$ est injective,
$\underline r$ est injective.\hfill\break
\hbox{\bf ii.} Si $[r]$ est bijective, $\overline r$ et $\underline r$ sont
bijectives. Donc :
$\overline r$ est une \'equivalence d'espaces de Hilbert. (En
particulier, deux pr\'esentations hilbertiennes de $H$ fournissent le
m\^eme noyau et la m\^eme r\'eduction de $H$).

\bigskip

\noi {\bf D\'emonstration} (de 6.2){\bf .} La premi\`ere assertion est
\'evidente. Pour \'etablir la seconde et montrer que 
$\overline r :
(Z/\overline{\delta C}) \to (Z'/\overline{\delta'C'})$ est injective,
if faut v\'erifier la propri\'et\'e suivante : si
$r(z) \in (\overline{\delta'C'})$, alors $z \in (\overline{\delta C})$.
Remarquons que, puisque $[r]$ est surjective, l'application 
$$-\delta' \oplus r : C' \oplus Z \to Z'$$ 
est surjective. Soit $K$ son noyau. Alors $(-\delta' \oplus r)$ admet un 
rel\`evement continu $\varphi: Z \to K^\bot$ qui est un isomorphisme
d'espaces de Hilbert (o\`u $K^\bot$ est l'orthogonal de $K$)
$$K = \{(c',\xi)\mid r(\xi) = \delta'(c')\}.$$
Soit alors $(z+\overline{\delta C}) \in \hbox{Ker}(\overline r)$ ; on a
donc : $r(z) = z' \in
(\overline{\delta'C'})$, et $z' = \lim \delta'(c'_n)$.
Soit $(\gamma'_n,z_n) := (c'_n,z) -
\varphi\big(z'-\delta'(c'_n)\big) \in K$, on a
$\varphi\big(z'-\delta'(c'_n)\big) \to 0$ et $z_n \to z$, en particulier.
Or, 
$$
-\delta'(\gamma'_n)+r(z_n)=(-\delta'\oplus r)(\gamma'_n,z_n)=
-\delta'(c'_n)+r(z)-\big(z'-\delta'(c'_n)\big)=r(z)-z'=0,
$$
donc $r(z_n)=\delta'\gamma'_n$, et $z_n \in \delta C$ par injectivit\'e
de $[r]$. On a donc bien $z \in (\overline{\delta C})$ comme 
annonc\'e.\hfill $\hbox{\cqfd}$

\bigskip

\noi {\bf 6.3. $\Gamma$-pr\'esentation.} Soit $\Gamma$ un groupe discret
agissant sur
l'espace vectoriel complexe $H$. Une $\Gamma$-pr\'esentation hilbertienne
$\delta : C \to
Z$ de $H$ est une pr\'esentation hilbertienne telle ue $\Gamma$ agisse de
mani\`ere unitaire
et \'equivariante sur $C$ et $Z$, l'action sur le quotient $(Z/\delta C)$
\'etant celle sur $H$.
Une telle pr\'esentation munit $\underline H$ et $\overline H$ d'une action
de $\Gamma$,
qui est unitaire sur $\overline H$.

\bigskip

\noi {\bf 6.4. $\Gamma$-dimension.} Soit $V$ un espace de Hilbert muni
d'une action unitaire de $\Gamma$. On dit que $V$ est de
$\Gamma$-dimension finie s'il existe un sous-ensemble fini
$\{v_1,\dots,v_m\}$ de $V$ tel que le sous-espace vectoriel engendr\'e
par les $(\gamma.v_i)$ $(\gamma \in \Gamma, 1 \le i \le m)$ soit dense
dans $V$. De mani\`ere \'equivalente, $V$ est un quotient ou un
sous-espace de $\big(\ell^2(\Gamma) \dis\otimes_{\Bbb C} {\Bbb
C}^m\big)$, ceci de mani\`ere compatible avec les actions naturelles
de $\Gamma$ (triviale sur ${\Bbb C}^m$).

\noi Dans cette situation, on d\'efinit la $\Gamma$-dimension
$\dim_\Gamma(V)$, qui est un nombre r\'eel $\big($inf\'erieur ou
\'egal \`a $m$, ici, et ind\'ependant du choix des $(v_i)\big)$. Voir
[P] pour cette notion. Les propri\'et\'es fondamentales (utilis\'ees
ici) de cette notion sont les suivantes :
\smallskip\noi
(6.4.1) \quad $\dim_\Gamma(V) \ge 0$, avec \'egalit\'e si et seulement 
si $V=0$.
\smallskip\noi 
(6.4.2) \quad Si $V$ est isomorphe \`a un sous-espace dense de $W$, alors
$$\dim_\Gamma V = \dim_\Gamma W.$$
(6.4.3) \quad Si $V = W \oplus W'$ (somme directe orthogonale), alors :
$$\dim_\Gamma V = \dim_\Gamma W + \dim_\Gamma W'.$$
(6.4.4) \quad $\dim_\Gamma\big(\ell^2(\Gamma)\big) = 1$.
\smallskip\noi
(6.4.5) \quad Si $\Gamma$ est un groupe fini de cardinal $|\Gamma|$, alors~ 
$\dim_\Gamma V=|\Gamma|^{-1} \dim_{\Bbb C}V$.

\bigskip

\noi {\bf 6.5. $\Gamma$-pr\'esentation finie.} Soit $\gamma : C \to Z$ une
$\Gamma$-pr\'esentation de l'espace vectoriel complexe $H$, muni d'une
$\Gamma$-action. On dit que $(\delta : C \to Z)$ est une
$\Gamma$-pr\'esentation finie de
$H$ si $C$ et $Z$ sont de $\Gamma$-dimensions finies. Alors
$(\overline{\delta C})$ est
de $\Gamma$-dimension finie au plus \'egale \`a celle de $C$, et 
$\overline H$ est aussi de
$\Gamma$-dimension finie. En fait, $\dim_\Gamma\overline H=\dim_\Gamma Z -
\dim_\Gamma\overline{\delta C}$. On pose alors 
$\dim_\Gamma H=\dim_\Gamma\overline H$.
\medskip
\noi Remarquons que cette d\'efinition peut \^etre donn\'ee m\^eme si
$C$ n'est pas suppos\'e \^etre de $\Gamma$-dimension finie.
\vfill\eject

\centerline{\bigbf Bibliographie}
\bigskip

\item{\bf [AG]} {\bf A.\ Andreotti, H.\ Grauert}, {\sl Th\'eor\`emes de 
finitude pour la cohomologie des espaces complexes}. Bull.\ Soc.\ Math.\ 
France {\bf 90} (1962), 193--259.
\medskip

\item{\bf [AV]} {\bf A.\ Andreotti, E.\ Vesentini}, {\sl Carleman estimates 
for the Laplace-Beltrami equation in complex manifolds}. Publ.\ Math.\ 
I.H.E.S.\ {\bf 25} (1965), 81--130.
\medskip

\item{\bf [A]} {\bf M.\ Atiyah}, {\sl Elliptic operators, discrete groups
and Von Neumann algebras.} Ast\'e\-risque 32-33 (1976), 43-72.
\medskip

\item{\bf [D]} {\bf J.-P.\ Demailly}, {\sl Estimations $L^2$ pour l'op\'erateur
$\overline
\partial$ d'un fibr\'e vectoriel hermitien semi-positif.} Ann. Sc. ENS 15
(1982), 457-511.
\medskip

\item{\bf [E]} {\bf P.\ Eyssidieux}, {\sl Th\'eorie de l'adjonction $L^2$ sur
le rev\^etement universel.} Preprint (1997).
\medskip

\item{\bf [G]} {\bf M.\ Gromov}, {\sl K\"ahler hyperbolicity and $L^2$-Hodge
theory.} J.\ Diff.\ Geom.\ {\bf 33} (1991), 263-291.
\medskip

\item{\bf [H]} {\bf L.\ H\"ormander}, {\sl An introduction to Complex 
Analysis in several variables}. 1st edition, Elsevier Science Pub., 
New York, 1966, 3rd revised edition, North-Holland Math.\ library, 
Vol 7, Amsterdam (1990).
\medskip

\item{\bf [K]} {\bf J.\ Koll\'ar}, {\sl Shafarevitch maps and automorphic
forms.} Princeton
University Press (1995).
\medskip

\item{\bf [Oh]} {\bf T.\ Ohsawa}, {\sl A reduction theorem for cohomology 
groups of very strongly $q$-convex K\"ahler manifolds}. Invent.\ Math.\ 
{\bf 63} (1981) 335-354$\,$; {\bf 66} (1982) 391-393.
\medskip

\item{\bf [P]} {\bf P.\ Pansu}, {\sl An introduction to $L^2$ Betti
numbers}. Preprint (1994).
\medskip

\item{\bf [W]} {\bf A.\ Weil}, {\sl Introduction \`a l'\'etude des vari\'et\'es
k\"ahl\'eriennes}.
Hermann (1958).
\vskip30pt
\noi
Fr\'ed\'eric Campana\hfil\break
Universit\'e de Nancy I, 
Facult\'e des Sciences, D\'epartement de Math\'ematiques\hfil\break
BP 239, 54506 Vandoeuvre les Nancy, France\hfil\break
E-mail: {\tt campana@iecn.u-nancy.fr}
\vskip10pt
\noi
Jean-Pierre Demailly\hfil\break
Universit\'e de Grenoble I, Institut Fourier, UMR 5582 du CNRS\hfil\break
BP74, 100 rue des Maths, 38402 Saint-Martin d'H\`eres, France\hfil\break
E-mail: {\tt demailly@ujf-grenoble.fr}

\end